%%Version 010901 %%
%%Corrections 101001%%

\input amstex
\documentstyle{amsppt}
\magnification=1200
\parskip 6pt
%\nologo
\NoBlackBoxes
%\TagsOnRight
\pagewidth{5.4in}
\hfuzz=5pt

%\macros
\def\sR{\hbox{I\kern-.1667em\hbox{R}}}
\def\R{\hbox{I\kern-.1667em\hbox{R}}}
\def\N{\hbox{I\kern-.1667em\hbox{N}}}
\def\Z{{\bold Z}}

\def\coog{C^0_0(G)}
\def\cog{C^0(G)}
\def\cwg{C^1(G)}
\def\cwog{C^1_0(G)}
\def\pair#1#2{\langle #1, #2 \rangle}

\def\la{\langle}
\def\ra{\rangle}
\def\lap{\Delta}
\def\did{\Delta_{iD}}
\def\o{{\Cal O}}
\def\dist{\text{\rm dist}}

\def\spec{\text{\rm spec}}

\def\mspec{\text{\rm mspec}}
\def\pspec{\text{\rm pspec}}
\def\deter{\text{\rm det}}

%new sjf definitions

%new rm definitions
\def\pair#1#2{\langle #1,#2 \rangle}

\def\sone#1#2{\left\{ \matrix 
                        #1 \\
                        #2  \endmatrix \right\}}
\def\stwo#1#2{\left[ \matrix 
                        #1 \\
                        #2  \endmatrix \right]}
\def\choose#1#2{\left( \matrix 
                        #1 \\
                        #2  \endmatrix \right)}
\def\chone{{\bold {1}}}
\def\ddx{\frac{d}{dx}}

\topmatter

\title
Diffusions on Graphs, Poisson Problems and Spectral Geometry
\endtitle

\author 
Patrick McDonald and Robert Meyers 
\endauthor

\affil
New College of Florida and The Courant Institute of Mathematical Sciences 
\endaffil

\address
New College of Florida and The Courant Institute of Mathematical Sciences 
\endaddress

\email  ptm$\@$virtu.sar.usf.edu,  meyersr$\@$cims.nyu.edu
\endemail

\abstract
We study diffusions, variational principles and associated
boundary value problems on directed graphs with natural weightings.  We
associate to certain subgraphs (domains) a pair of sequences, each of
which is invariant under the action of the automorphism group of the
underlying graph.  We prove that these invariants differ by an
explicit combinatorial factor given by Stirling numbers of the first
and second kind.  We prove that for any domain with a natural
weighting, these invariants determine the eigenvalues of the Laplace
operator corresponding to eigenvectors with nonzero mean.  As a 
specific example, we investigate the relationship between our
invariants and heat content asymptotics, expressing both as special
values of an analog of a spectral zeta function.

\endabstract

\date June 28, 2001\enddate

\endtopmatter

\document

\heading
1:  Introduction
\endheading

In this note we study relationships between combinatorial analogs of
various probabilistic and geometric objects associated to domains in a
Riemannian manifold.  To precisely state our results, we begin by
establishing some notation.

Let $G= (V,E)$ be a connected oriented bidirected graph with vertex
set $V$ and edge set $E.$   Given $e \in E,$ we will represent $e$ as
an ordered pair $e = (t(e), h(e))$ where $t(e), \ h(e) \in V$ (see
section 2 for definitions).   

We will denote by $\cog$ the vector space of real valued functions on
$V$ and by $\cwg$ the vector space of real valued functions on $E.$
There is a natural coboundary operator $d:\cog \to \cwg$ defined
by $$d f(e) =f(h(e)) - f(t(e)).$$
Let $\coog\subset \cog$ and $\cwog \subset \cwg$ be the subspaces
consisting of those functions with compact support. 

Let $W_V: V \to \sR^+$ be a vertex weighting.  Associated to $W_V$
there is an inner product on $ \coog$ defined by 
$$\align
\la f,g\ra_V & = \sum_{x\in V} f(x) g(x) W_V(x). \tag{1.1}
\endalign$$
Similarly, a function $W_E: E \to \sR^+$ defines an inner product
$\la \cdot \ , \ \cdot \ra_E$ on $\cwog.$  

Given a pair of functions $W_V$ and $W_E$ as above, we call the
ordered pair $W = (W_V,W_E)$ a weighting for $G$ if 
$$\align
W_E(x,y)W_V(x) & = W_E(y,x)W_V(y) \tag{1.2}
\endalign$$
where $x, y  \in V.$  If $\o$ is an orientation for $G$ (cf section 2)
and $W$ is a weighting, we call the triple $(G,\o,W)$ a {\it graph with
geometry}.  

Given a weighting $W,$ we will denote by $d^*_W:\cwg \to \cog$ the
associated adjoint of the coboundary map $d.$  We will denote by
$\Delta : \cog \to \cog$ the (vertex) Laplacian, $\Delta = d^*_W d$
(we will not investigate the edge Laplacian $d d^*_W$).  

Interest in applications involving discrete Laplace operators of the
type defined above can be traced to Kirchoff \cite{K} who modeled simple
circuits as finite graphs with each edge corresponding to the
conductance of a given circuit component (cf \cite{DSn} for a survey
of random walks and electrical networks).  As Kirchoff established, it
is possible to give graph theoretic formulations for the laws
governing current behavior in a simple circuit (Kirchoff's laws for
voltage and current, Ohm's law). Given a simple circuit (ie a graph
with each edge weighted to represent the conductance of a given
component and each vertex of weight one) and a given input current, it
is possible to formulate the problem of finding the induced circuit
current  as a Dirichlet problem involving the edge Laplacian.  In
addition, there is a solution to this Dirichlet problem given by
``energy minimization'' (Thomson's Principle, cf \cite{B}).   

Since Kirchoff, the study of graph Laplacians has yielded a remarkable
wealth of information in a variety of contexts (the references
\cite{B}, \cite{C} and references therein provide expository
introductions to some of these applications).  Among those fields 
where graphs and their associated discrete boundary value problems
have found interesting applications are potential theory (cf \cite{B},
\cite{Du}, and references therein), spectral theory (cf \cite{DS},
\cite{Ge}, \cite{C}), and differential geometry and global analysis (cf
\cite{Do}, \cite{F1}, \cite{F2},  \cite{V1}, \cite{V2}).  In this
paper, we focus on several such applications.  To be more precise,
let $W=(W_V,W_E)$ be a weighting for $G$ and define the associated
vertex weighting as $$w_V(x) = \sum_{(x,y)\in E} W_E(x,y).$$  Let 
$$\align
p(x,y) & = \frac{W_E(x,y)}{w_V(x)} \tag{1.3}
\endalign$$
be the transition probabilities for a random walk, denoted $X_n,$ on
the vertices of $G.$  For $x\in V,$ let 
$P^x$ be the associated measure charging trajectories beginning at $x$
and let $E^x$ be the corresponding expectation operator (a classical
reference for random walks is \cite{Sp}; \cite{A} is a survey
concerning the role of random walks in a variety of contexts).  Let
$D\subset G$ be a domain of $G$ (cf Definition 2.2) and let $\eta$ be
the index of the first exit time from the interior of $D$ (cf
Definition 2.4): 
$$\align
\eta & = \inf \{n\geq 0: X_n \notin D\}. 
\endalign$$
Rescaling transition times to reflect graph geometry, let $\tau,$ the
first exit time from the interior of $D,$ be defined by (cf   
Definition 2.4): 
$$\align
\tau & = \sum_{n=0}^{\eta -1} \frac{1}{w_V(X_n)}. \tag{1.4}
\endalign$$
Many of the applications of the study of graph Laplacians to 
geometry and analysis involve the relationship of exit times and
hitting times for random walks on $G$ to the discrete Laplace operator
and a corresponding theory of boundary value problems.  Our first
result is a consequence of identifying the precise relationship
between exit time moments of the natural random walk described above
and solutions to natural Poisson problems involving graph Laplacians
on $G$ (cf Theorem 4.1).  More precisely, given a domain $D$ of $G$
and a positive integer $k$, define real numbers $A_{1,k}= A_{1,k}(D)$
and $A_{2,k}=A_{2,k}(D)$ by  
$$\align
A_{1,k} & = \pair{E^x[\tau^k]}{1}_V \tag{1.5} \\
A_{2,k} & = k!\sup_{f\in C^0(iD)\setminus \{0\}}
\frac{ \pair{f}{1}_{iD}^2}{ |\pair{f}{\Delta_{iD}^kf}_{iD}|}
\tag{1.6}  
\endalign$$
where $iD$ is the collection of interior vertices of $D$ (cf
Definition 2.2), $\did$ is the interior Laplace operator (cf
\rom{(2.9)} and Proposition 2.3) and the inner product in \rom{(1.6)}
is induced by restriction to functions supported on $iD$ (cf section
2).  For reasons that will soon become apparent, we will refer to the
sequence $\{A_{1,k}(D)\}$ as the {\it moment spectrum} of $D$ and we
will write $\hbox{mspec}(D) = \{A_{1,k}(D)\}.$  Similarly, we will
refer to  the sequence $\{A_{2,k}(D)\}$ as the {\it Poisson spectrum}
of $D$ and we will write $\hbox{pspec}(D) = \{A_{2,k}(D)\}.$ We prove
the following  

\proclaim{Theorem 1.1}Suppose that $(G,\o,W)$ is a connected graph
with geometry (cf Definition \rom{2.1}).  Let $p:V\times V 
\to \R$ be the transition probabilities associated to the weighting
$W$ as defined in \rom{(1.3)} above, and let $X_n$ be the random walk
determined by the transition probabilities $p(x,y).$  Let $D$ be an
$\alpha$-weight regular domain (cf Definition \rom{(2.5)}).  For $k$ a
positive integer, let $A_{1,k}$ and $A_{2,k}$ be 
defined by \rom{(1.5)} and \rom{(1.6)}, respectively.  Then 
$$\align
A_{1,k} & = A_{2,k} + \sum_{j=1}^{k-1}
(-\alpha)^{j-k}\sone{k}{j} A_{2,j}  \tag{1.7}
\endalign$$
where the coefficients $\sone{k}{j}$ are the Stirling numbers of
the second kind (cf Definition \rom{(4.1)}).  Similarly, 
$$\align
A_{2,k} & = A_{1,k} + \sum_{j=1}^{k-1}
(-\alpha)^{j-k}\stwo{k}{j} A_{1,j}  
\tag{1.8} 
\endalign$$
where the coefficients $\stwo{k}{j}$ are the Stirling numbers of the
first kind (cf  Definition \rom{(4.1)}). 
\endproclaim

Theorem 1.1 implies that for every weight regular domain $D,$
$\mspec(D)$ determines $\pspec(D)$ and vice versa. 

Theorem 1.1 is motivated in part by an analogous result concerning
Brownian motion on smoothly bounded domains in Riemannian manifolds.
The precise statement is as follows:  Let $M$ be a manifold with
Riemannian metric $g.$  Suppose that $X_t$ is classical Brownian
motion in $M,$ that $D\subset M$ is a smoothly bounded domain with
compact closure, and that $\tau$ is the first exit time of Brownian
motion from $D.$  Let $P^x$ be the measure charging paths beginning at
$x\in M$ and for $k$ a positive integer, let $B_{1,k}$ be defined by 
$$\align
B_{1,k} & = \int_D E^x[\tau^k] dg \tag{1.9}
\endalign$$
where $E^x$ is expectation with respect to $P^x$ and $dg$ is the
metric density.  As a function of the domain, the sequence
$\{B_{1,k}(D)\}$ has a number of interesting properties.  For domains
in the plane, $B_{1,1}(D)$ is known as the {\it torsional rigidity}
associated to a beam of uniform cross-section $D$ and arises in the
theory of elasticity.  In general, the sequence $\{B_{1,k}(D)\}$ is 
invariant under the action of the isometry group of $M$ and hence
contains geometric information for the domain $D$ (see \cite{KMM},
\cite{M} for related work).  Define  
$$\align
B_{2,k} & = k! \sup_{f \in {\Cal F}_k} \frac{\left(\int_D f
dg\right)^2}{\left| \int_D f {\Cal L}^k f dg\right|}\tag{1.10}
\endalign$$
where ${\Cal L}$ is one half the Laplace operator and ${\Cal F}_k =
\{f\in C^\infty(D): \lap^j f = 0 \hbox{ on } \partial D, \ \ 0 \leq j
<k\}.$  Then for all positive integers $k$ (cf \cite{KM} for the
Euclidean case, \cite{M} for the general case)  
$$\align
B_{1,k} & = B_{2,k}. \tag{1.11}
\endalign$$
Suppose we are given a Riemannian manifold $M.$  A
triangulation of $M$ leads to a bidirected graph with orientation, as
well as approximations of the associated Laplace operator (by finite
difference operators) and the Brownian motion $X_t$ (by random walks).
These approximations give rise to analogs ($A_{1,k}$) of the geometric
invariants $B_{1,k}$ and similarly for the variational quotients.
Theorem \rom{1.1} precisely quantifies the difference which arises in
approximating two geometric invariants in terms of the Stirling cycle
numbers, combinatorial objects which arise naturally in a variety of
contexts in graph theory and statistical mechanics. 

The remainder of our results demonstrate that the invariants $A_{1,k}$
and $A_{2,k}$ play an important role in the analysis 
and geometry of infinite graphs.  These results provide a
relationship between the spectrum of the interior Laplace operator
associated to a domain, the Poisson spectrum of the domain, and the
moment spectrum of the domain.  We begin with a definition:

\proclaim{Definition 1.1} Let $(G,\o,W)$ be a graph with geometry
and let $D$ be a domain in $G.$  Let $\spec(D)$ be the spectrum of
the interior Laplace operator associated to $D.$  We define the set
$\spec^*(D)$ by 
$$\align
\spec^*(D) & = \{\lambda \in \spec(D):
\pair{\phi_\lambda}{\chone_{iD}}_V \neq 0 \} \tag{1.12}
\endalign$$
where $\phi_\lambda$ is a normalized eigenvector associated to
$\lambda$ and $\chone_{iD}$ is the indicator function of $iD.$
\endproclaim
We emphasize that $\spec^*(D)$ contains no information concerning
spectral multiplicity; it is a subset of the real numbers consisting
of those eigenvalues whose corresponding eigenspace projects
nontrivially onto constant functions.

Our second result can now be concisely stated:

\proclaim{Theorem 1.2}Suppose that $(G,\o,W)$ is a graph with geometry
and that $D$ and $D'$ are domains in $G.$  Then, with the notation as
in Definition \rom{1.1}, 
$$\align
\pspec(D) = \pspec(D') & \hbox{ implies } \spec^*(D) =
\spec^*(D') 
\endalign$$
and we say that $\pspec(D)$ determines $\spec^*(D).$  Moreover,
suppose ${\Cal D}_N$ is the collection of all domains for which the
cardinality of $\spec^*(D)$ is $N:$
$$\align
{\Cal D}_N & = \{D \subset G: D \hbox{ a domain }, \
\text{\rm card}(\spec^*(D)) = N\}. \tag{1.13}
\endalign$$
Then there are $N$ rational functions $f_i, \ 1 \leq i \leq N, \
f_i:\R^{2N} \to \R$ such that for every $D \in {\Cal D}_N,$ the
roots of the polynomial 
$$P_N(x) = x^N + \sum_{i=0}^{N-1} f_i(A_{2,0}(D), A_{2,1}(D), \dots,
A_{2,2N-1}(D))x^i \tag{1.14} $$
give the elements of $\spec^*(D).$  When the domains under
consideration are $\alpha$-weight regular, the same claims holds with
$\mspec(D)$  replacing $\pspec(D).$  
\endproclaim

Thus, given a domain $D,$ any invariant defined in terms of
$\spec^*(D)$ can be defined in terms of $\pspec(D)$ (and/or
$\mspec(D)$).  An important example of such invariants are given by
the notion of {\it heat content}.  Recall, if $D$ is a smoothly
bounded domain in a Riemannian manifold $M$ and $H(x,t)$ solves the
heat equation with intial data
$$\aligned
\Delta H & = \partial_t H \hbox{ on } D\times (0,\infty) \\
H(x,0) & = 1 \hbox{ on } D \\
H(y,t) & = 0 \hbox{ on } \partial D \times (0,\infty)
\endaligned\tag{1.15}$$
then the heat content of $D$ is the function $Q(t)$ defined by 
$$\align
Q(t) & = \int_D H(x,t) dg \tag{1.16}
\endalign$$
where $dg$ is the metric density.  By a theorem of Gilkey \cite{Gi},
$Q(t)$ admits an asymptotic expansion for small $t:$
$$\align
Q(t) & \simeq \sum_{n=0}^\infty q_n t^n \tag{1.17}
\endalign$$
where the coefficients $q_n$ are given as integrals of metric
invariants associated to $D.$  The coefficients in \rom{(1.17)} are
sometimes referred to as the {\it heat content asymptotics}.  We prove 

\proclaim{Theorem 1.3} Suppose that $(G,\o,W)$ is a graph with
geometry and that $D$ is a domain of $G.$  Suppose the heat content
asymptotics of $D$ are given by $\{q_n\}$ (cf Definition
\rom{(5.2)}).  Then $\{q_n\}$ determines $\spec^*(D).$  In addition,
$\{q_n\}$ is determined by $\pspec(D),$ and if $D$ is $\alpha$-weight
regular, the same is true with $\pspec(D)$ replaced by $\mspec(D).$    
\endproclaim

We note that there are natural questions arising from our results.
For example,  to what extent can we establish converse statements for
the results of Theorem 1.2 and Theorem 1.3 (ie, to what extent
does $\spec^*(D)$ determine $\pspec(D)$)?  We address these issues in
Corollary 5.1 and Corollary 5.3, respectively, where we prove that
$\spec^*(D)$ together with spectral information related to the volume
of $D$ (cf Definition 5.1) determine the Poisson spectrum, and
similarly for the heat content asymptotics.  In Definition
5.3 we define a variant of the spectral zeta function for the
Dirichlet Laplacian.  Our final result, Corollary 5.4, establishes
that special values of our zeta function give both the invariants
defined in \rom{(1.6)} (values at positive integers) and the heat
content asymptotics (values at negative integers).

To prove our theorems, we study a natural family of Poisson problems on
the domain $D$ and develop the potential theory required to analyze
such systems.  Of particular importance is the precise relationship
between the geometry of the underlying graph, the transition operator,
the inverse of the Laplace operator, and exit time moments for the
natural random walk on the graph (cf Corollary \rom{(2.1)}, Lemma
\rom{(4.2)} and Theorem \rom{(4.1)}).  To prove Theorem 1.1 we prove
that the variational quotients appearing in \rom{(1.6)} are given by
integrating the solutions of the Poisson problems arising in our study
(cf Theorem 3.2), we establish an analog of Dynkin's formula
(cf Theorem 3.1) and we use this to relate the moments appearing in
\rom{(1.5)} to the solutions of our Poisson problems (cf Theorem 4.1).
To prove Theorems 1.2 and 1.3 we rely on the solution of the Stieltjes
Moment Problem and explicit computations.  Our computations provide
algebraic descriptions of the invariants of interest and indicate the
precise nature of the relationship between the sets of invariants.   

The paper is organized as follows:  In the second section we review 
the necessary machinery of probability, analysis and graph theory.  In
the third section we examine iterated solutions of discrete Poisson
problems.  In the fourth section we present the proof of Theorem 1.1.
In the fifth and final section we present the proofs of Theorem 1.2
and Theorem 1.3, as well as a number of corollaries involving converse
results.

\heading
2: Background and Notation
\endheading

Let $G= (V,E)$ be a connected bidirected graph with vertex set $V$ and
edge set $E\subset V\times V.$   Given $e \in E,$ we will represent
$e$ as an ordered pair $e = (t(e), h(e))$ where $t(e), \ h(e) \in V$
($h(e)$ is the terminal point of $e$ and $t(e)$ is the initial point
of $e$).  We will say that a vertex $x$ is incident with an edge $e$
if $x=t(e)$ or $x=h(e).$  For $x$ a vertex,  we denote by $\deg(x)$
the degree of $x:$  $\deg(x) = \hbox{number of } e$ such that $x$ is
incident with $e.$ We will restrict our attention to graphs which
admit no self edges (ie, edges of the form $(x,x).$)  We will denote
the statement ``$x$ is adjacent to $y$'' by $x \sim y.$   

Given an edge $e = (x,y),$ we will denote the {\it opposite edge}
$(y,x)$ by $e^{-1}.$  An {\it orientation} of a bidirected graph $G$
is a subset $\o \subset E$ satisfying $e \in \o$ if and only if
$e^{-1} \notin \o.$  

A vertex weighting (or volume) on a bidirected graph $G=(V,E)$ is a
function $W_V : V \to \R^+.$  An edge weighting on a bidirected graph
$G= (V,E)$ is function $W_E :V \times V \to \R^{\geq 0}$ which is
supported on $E.$  

\proclaim{Definition 2.1} Suppose that $G$ is a connected bidirected graph
with orientation $\o,$ vertex weight $W_V,$ and edge weight $W_E.$  We
say that the pair $W = (W_V,W_E)$ is a weighting for $G$ if 
$$\align
W_E(x,y)W_V(x) & = W_E(y,x)W_V(y) \tag{2.1}
\endalign$$
for all $x,y \in V.$  Given a weighting $W,$ we associate to $G$ a
second vertex weighting, $w_V(x),$ called the natural auxillary
weighting, defined by 
$$\align
w_V(x) & = \sum_{(x,y) \in E} W_E(x,y). \tag{2.2}
\endalign$$
A triple $(G,\o,W)$ where $G$ is a connected bidirected graph, $\o$ is
an orientation of $G$ and $W$ is a weighting for $G$ is called a {\it
graph with geometry}.  
\endproclaim

As in the introduction, we will denote by $\cog$ the vector space of
real valued functions on $V$ and by $\cwg$ the vector space of real
valued functions on $E.$ There is a natural coboundary operator
$d:\cog \to \cwg$ defined by $$d f(e) =f(h(e)) - f(t(e)).$$Let
$\coog\subset \cog$ and $\cwog \subset \cwg$ be the subspaces
consisting of those functions with compact support.  

Let $(G,\o,W)$ be a graph with geometry.  The weighting $W$ gives rise
to a pair of inner products 
$$\align
 & \pair{\cdot}{\cdot}_V : \coog \times \coog \to \R \\
 & \pair{\cdot}{\cdot}_E : \cwog \times \cwog \to \R
\endalign$$
defined by
$$\align
\la f,g\ra_V & = \sum_{x\in V} f(x) g(x)W_V(x) \tag{2.3}\\
\la F,G\ra_E & = \sum_{x\in V} \sum_{(x,y)\in \o} F(x,y) G(x,y)
W_E(x,y) W_V(x). \tag{2.4}
\endalign$$
The inner products associated to $W$ give rise to an adjoint map
$$d_W^*: \cwog \to \coog$$ which in turn gives rise to the (vertex)
Laplace operator $\Delta : \coog \to \coog$ defined by 
$$\align 
\Delta & = d_W^* d. 
\endalign$$

\demo{Example 2.1}  Let $V=\Z^n$ be the integer lattice
in $\R^n$ (with the standard Riemannian structure).  Define an edge
structure by setting 
$$\align
E & = \{(x,y): x, \ y \in V, \ \dist(x,y) \leq 1\}
\endalign$$
where $\dist(x,y)$ is the standard Euclidean distance between $x$ and
$y.$  The orientation of $\R^n$ induces an orientation on the graph
$G=(V,E):$ 
$$\align 
 & (x,y) \in \o \hbox{ if the tangent to the line from $x$ to $y$ is
positively oriented. }
\endalign$$
Define an edge weighting by
$$\align
W_E(x,y) & = \cases
             \frac{1}{(\dist(x,y))^2} & \hbox{ if } (x,y) \in E \\
             0 & \hbox{ elsewhere }
             \endcases 
\endalign$$
Define a vertex weighting by 
$$\align
W_V(x) & = 1.
\endalign$$
Then $W = (W_V,W_E)$ is a natural weighting.  The corresponding
Laplace operator is the standard discretization of the Laplace
operator acting on functions on $\R^n.$

To generalize the above, suppose that $M$ is an oriented
complete Riemannian manifold with metric $g$ and injectivity radius
bounded below by a positive constant, $b.$  Choose a discrete
collection of points $\{x_i\} =V \subset M$ satisfying
\roster
\item $\dist_g(x_i,x_j) \geq \frac{b}{4}$ for $i\neq j$
\item For every $x\in M$ there is some $i$ such that $\dist_g(x,x_i) <
\frac{b}{2}$ 
\endroster
where $\dist_g$ is the $g$-induced distance function on $M.$  Define
an edge structure on $V$ by the relationship 
$$\align
x_i \sim x_j  &  \hbox{ if } \dist_g(x_i,x_j) < b \hbox{ and } i
\neq j.
\endalign$$
Then $G=(V,E)$ is a connected bidirected graph.  The orientation on
$M$ induces an orientation on $G:$  Given an edge $(x,y)$ we say $(x,y)
\in \o$ if the tangent vector at $x$ to a distance minimizing geodesic
connecting $x$ to $y$ is given by the orientation on $M.$  We define
an edge weighting on $E$ and a vertex weighting on $V$ by 
$$\aligned
W_E(x,y) & = \frac{1}{(\dist_g(x,y))^2} \\
W_V(x) & = 1.
\endaligned$$
From the definition, it is clear that $W=(W_V,W_E)$ is a weighting
and that the triple $(G,\o,W)$ is a graph with geometry.  In the case
of Euclidean space $\R^n,$ taking $V = Z^n$ where $\Z^n$ is the
integer lattice and $b = \sqrt{2},$ results in the triple being the 
standard nearest-neighbor graph on the given lattice where each vertex
and each edge is assigned the same weight.  The associated Laplace
operator is the standard finite difference operator.

In \cite{V1} Varopolous constructs a similar transition operator and
proves that as the mesh described above becomes fine, the natural
random walk converges to classical Brownian motion on $M,$ and uses
this fact to relate transience of Brownian motion on regular coverings
of a compact Riemannian manifold to transient properties of the
corresponding group of deck transformations.  In addition, Varopolous
uses the random walk to give large time bounds on the heat kernel of a
noncompact manifold.  
\enddemo

Our first proposition characterizes the action of the Laplacian as a
weighted average (cf \cite{Do} for the case in which each vertex has
the same degree and the edge weights are constant)

\proclaim{Proposition 2.1} Suppose that $(G,\o,W)$ is a graph with
geometry, $\Delta$ the associated Laplace operator. The action of 
$\lap$ on a function $f$ at a vertex $x$ is given by
$$-\lap f(x) = \sum_{y \in V} f(y)W_E(x,y) - w_V(x) f(x)
\tag{2.5} $$ 
where $W_E$ is the edge weight function and $w_V$ is the auxillary vertex
weight function \rom{(2.2)}
\endproclaim

\demo{Proof}  For $x \in V,$ let $\delta_x \in \coog$ be defined by 
$$\aligned
\delta_x(y) & = \cases 1 & y=x \\
                       0 & \hbox{elsewhere}  \endcases
\endaligned\tag{2.6}$$
From the definition of the inner product \rom{(2.3)},
$$\align
\pair{\Delta f}{\delta_x}_V & = \Delta f(x) W_V(x). \tag{2.7}
\endalign$$
As $d^*_W$ is the adjoint of the coboundary operator with respect to
$W,$ we have, via \rom{(2.4)},
$$\align
\pair{\Delta f}{\delta_x}_V & = \pair{df}{d\delta_x}_E \\
 & = \sum_{(x,y)\in \o} (f(y)-f(x))(-1) W_E(x,y)W_V(x)  \\
 & \ \ + \sum_{(y,x)\in \o} (f(x) - f(y))W_E(y,x)W_V(y) .
\endalign$$ 
Using the fact that $\o$ is an orientation and relationship
\rom{(2.1)}, we obtain 
$$\align
\pair{\Delta f}{\delta_x}_V & =  \sum_{(x,y)\in E}
(f(x)-f(y))W_E(x,y)W_V(x) .
\endalign$$
Using the definition of the auxillary vertex weighting \rom{(2.2)} and
\rom{(2.7)}, the proof of the proposition follows.   
\enddemo

We next develop the machinery of boundary value problems for graphs
and subgraphs.

Suppose that $G= (V,E)$ is a graph, $(G,\o,W)$ a graph with
geometry.  If $A\subset V,$ we define 
$$\align
C^0(G,A) &  = \{f\in C^0(G): f|_{V \setminus A} = 0\}.
\endalign$$
There are natural isomorphisms $J_A:C^0(A) \to C^0(G,A)$ which induce
inclusions: 
$$\align
I_A & : C^0(A) \to C^0(G) 
\endalign$$
and projections:
$$\align
P_A & : C^0(G) \to C^0(A) .
\endalign$$
The vertex weighting $W_V$ induces a weighting on $A$ and an inner
product on compactly supported functions: 
$$\align
\pair{f}{g}_A & =  \pair{I_Af}{I_Ag}_V.
\endalign$$

\proclaim{Definition 2.2}  Suppose that $G= (V,E)$ is a bidirected graph.
Suppose that $D= (V',E')$ is a bidirected subgraph of $G.$  We say
that $x\in V'$ is an interior vertex of $D$ if, for all $y \in V$ with 
$(x,y)$ an edge of $G,$ we have that $y \in V'$ and $(x,y) \in E'.$  We denote
the collection of all interior vertices of $D$ by $iD.$  We call all
vertices of $D$ which are not interior vertices of $D$ {\it boundary
vertices} of $D.$  We denote the collection of all boundary vertices
of $D$ by $\partial D.$ A {\it domain} of $G$ is a finite connected
bidirected subgraph of $G$ with nonempty interior vertex set.  
\endproclaim

Note that if $(G,\o, W)$ is a graph with geometry and $D$ is a
domain of $G,$ then $\o$ induces (by restriction) an orientation,
$\o',$ of $D,$ and $W$ induces (by restriction) a weighting, $W',$ for
$D$ which coincides with $W$ at all interior vertices.  The triple
$(D,\o',W')$ gives rise to a Laplace operator, $\Delta_D: C^0(D) \to
C^0(D),$ the induced Laplace operator associated to $D.$  

\proclaim{Proposition 2.2} Suppose that $(G,\o,W)$ is a graph with 
geometry, $\Delta$ the associated Laplace operator.  Let $D$ be a
domain in $G,$ $\Delta_D$ the induced Laplace operator associated to
$D.$  Then for all $f \in C^0(D),$ for every interior vertex $x$ of
$D,$ 
$$\align
\Delta_D f(x) & = \Delta I_D f(x) \tag{2.8}
\endalign$$
where $I_D$ is the natural inclusion, $I_D: C^0(D) \to C^0(G).$ 
\endproclaim

\demo{Proof} If $D= (V',E'),$ and $(D,\o',W')$ denotes the domain with
induced structure, then following the proof of Proposition
\rom{(2.1)}, we have 
$$\align
\pair{\Delta_D f}{\delta_x}_{V'} & = \sum_{(x,y) \in E'} (I_Df(x) - I_Df(y))
W_{E'}(x,y)W_{V'}(x) .
\endalign$$
Since $W_V(x) = W'_V(x)$ and $W_E(x,y) = W'_E(x,y)$ at all interior
vertices, \rom{(2.8)} follows. 
\enddemo

Note that it is not necessarily the case that \rom{(2.8)} holds at
boundary vertices: the respective auxillary weightings at a boundary
vertex may not coincide.  Nonetheless, we can give a concise
description of the action of the domain Laplacian as follows:
enumerate the vertices of $D$ with the first $m$ vertices interior,
the remaining $N-m$ vertices boundary. For $x \in D,$ let $\delta_x$
be as given in \rom{(2.6)}. Then $\{\delta_{x_i}\}_{i=1}^N$ is a basis
for $C^0(D)$ and the domain Laplacian with respect to this basis has
the form: 
$$\align
\Delta_D & = \left( \matrix 
                      \Delta_{iD,iD} &  \Delta_{\partial D,iD} \\
                      \Delta_{iD,\partial D} &  \Delta_{\partial
D,\partial D} 
                 \endmatrix \right)
\endalign$$
where the action of the operators $\Delta_{A,B}$ is given by 
$$\align
\Delta_{A,B} & = P_B \Delta_D I_A \tag{2.9}
\endalign$$
where $I_A$ and $P_B$ are the natural inclusion and projection maps.
Of particular interest is the interior Laplace operator acting on
functions on $iD,$ $\Delta_{iD,iD}:
C^0(iD) \to C^0(iD).$  For notational ease, we will write
$\Delta_{iD,iD} = \did.$ 

\proclaim{Proposition 2.3} Suppose that $(G,\o, W)$ is a graph with
geometry and that $D$ is a domain of $G.$  Suppose that $\did$
is the interior Laplace operator defined by \rom{(2.9)} with $A=B=iD.$
Then, for all $f\in C^0(iD),$ for all $x \in iD,$
$$\align
-\did f(x) & = \sum_{y\in iD} f(y) W_E(x,y) - w_V(x) f(x) \tag{2.10}
\endalign$$
and the interior Laplace operator is self-adjoint.
\endproclaim

\demo{Proof} Let $x \in iD.$  By definition, $\did f(x) = P_{iD}
\Delta_D I_{iD} f(x).$ 
By Proposition 2.2 and Proposition 2.1,
$$\align
-\Delta_D I_{iD} f(x)  & = -\Delta I_D f(x) \\
 & = \sum_{y\in V} I_Df(y) W_E(x,y) - w_V(x) I_Df(x) .
\endalign$$
Since $I_Df$ is supported on $iD,$ \rom{(2.10)} follows.  That $\did$
is self-adjoint on $C^0(iD)$ follows immediately from \rom{(2.10)}.
\enddemo

To study properties of the Laplace operators defined above, we
develop two natural random walks.  It is expression \rom{(2.5)} and
the corresponding variants which link the behavior of the random
walks to solutions of boundary value problems involving the
associated Laplace operators.

\proclaim{Definition 2.3} Let $(G,\o, W)$ be a graph with 
geometry.  The transition probabilities associated to $(G,\o,W)$ is
the function $p:V\times V \to \R$ defined by
$$\align
p(x,y) & = \frac{W_E(x,y)}{w_V(x)}. \tag{2.11}
\endalign$$
\endproclaim

Note that for every vertex $x,$ the ratio $p(x,y)
=\frac{W_E(x,y)}{w_V(x)}$ is a probability distribution on the 
collection of vertices adjacent to $x.$  These transition
probabilities induce a natural random walk, $X_n,$ taking values in the
vertices of $G.$  Let $P^x, \ x \in V,$ be the associated family of
probability measures charging paths beginning at $x$ and let $E^x$
denote the corresponding expectations.   

Because $G$ is connected and $W_E$ is supported on $E,$ the walk $X_n$
is irreducible.  Because $W$ is a weighting, relationship \rom{(2.1)}
implies that the natural random walk is reversible with stationary measure
$\pi(x) = W_V(x)w_V(x).$  For the graph arising in Example 2.1, the
natural random walk is the usual simple symmetric random on the
vertices of the integer lattice in $\R^n.$

Associated to $X_n$ there is a transition operator $T:\coog \to \coog$
defined by 
$$\align
Tf(x) &  = \sum_{y\in V} p(x,y)f(y) \tag{2.12}
\endalign$$
Let $p_n(x,y)$ be defined by recursion:
$$\align
p_1(x,y) & = p(x,y) \\
p_n(x,y) & = \sum_{z\in V} p(x,z)p_{n-1}(z,y).
\endalign$$
Iterating $T$ and using Chapman-Kolmogorov gives 
$$\aligned
T^nf(x) & =  \sum_{y\in V} p_n(x,y) f(y).
\endaligned \tag{2.13}$$

\proclaim{Proposition 2.4} Let $(G,\o,W)$ be a graph with geometry,
$(X_n,P^x)$ the associated random walk, and $T$ the 
transition operator associated to $X_n.$  Let $\Omega_w: \coog \to
\coog$ be the multiplication operator defined by $\Omega_wf(x) =
w_V(x)f(x).$  Then
$$\align
\Delta & = \Omega_w(I-T) \tag{2.14}
\endalign$$
where $I$ is the identity operator.
\endproclaim
\demo{Proof}  This follows directly from the definition of $p(x,y)$ in
\rom{(2.11)}, the definition of $T$ in \rom{(2.12)}, and \rom{(2.5)}.  
\enddemo

In the case of manifolds, there is a close relationship between the
behavior of a Brownian process in a given domain contained in a
Riemannian manifold, and the solutions of boundary value problems for
the Laplace operator on the given domain.  In the sequel, we present 
the machinery required to formulate analogs of these notions in the
category of graphs with geometry.  The material is by no means new and
is included for completeness (cf \cite{C}, \cite{B}, \cite{Do}, and 
references therein). 

\proclaim{Definition 2.4} Let $(G,\o,W)$ be a graph with natural
geometry, $D$ a domain of $G$ with $\partial D \neq \emptyset.$  Let
$X_n$ be the natural random walk associated to $(G,\o, W).$  The
index of the first exit time of $X_n$ from $iD$ is defined by
$$\align
\eta & = \inf\{ n : X_n \notin iD\} \tag{2.15}
\endalign$$
The first exit time of $X_n$ from $iD$ is defined by
$$\align
\tau & = \sum_{n=0}^{\eta -1} \frac{1}{w_V(X_n)}. \tag{2.16}
\endalign$$
\endproclaim

Note that $\eta$ and $\tau$ are measurable with respect to the
filtration associated to $X_n.$  Our definition of $\tau$ is motivated
by our desire to model approximation procedures for the case of
manifolds.  We view $w_V(x)$ as a measure of the ``size'' of the
neighborhood defined by the vertex $x.$  Thus our choice reflects our
desire to correlate ``transition time for leaving a neighborhood
defined by a vertex $x$'' with ``the size of a neighborhood defined by
a vertex $x$.''

Note that if for all $x\in iD,$ $w_V(x) = \alpha, \ \alpha$ a
constant, then $\alpha \tau =  \eta$ and the walk transitions at a
constant rate, proportional to the index.

\proclaim{Definition 2.5} Let $(G,\o,W)$ be a graph with geometry, $D$
a domain of $G.$  We say that $D$ is $\alpha$-weight regular if for
all $x \in iD, \ w_V(x) = \alpha.$
\endproclaim

Let $(G,\o,W)$ be a graph with geometry, $D$ a domain of $G$
with $\partial D \neq \emptyset.$  Because $D$ is finite, $$P^x(X_n
\in iD \ \ \forall n) = 0.$$In 
addition, for all natural numbers $k,$ $$E^x[\tau^k] < \infty.$$ 

\proclaim{Definition 2.6} Let $(G,\o,W)$ be a graph with geometry, $D$
a domain of $G$ with nonempty boundary.  Let 
$X_n$ be the process associated to $(G,\o,W).$  Let $x^*$ be a
cemetary point not contained in $D.$  The process $X_n$
killed upon leaving $iD$ is the process with transition 
probabilities $q(x,y)$ given by 
$$\aligned
q(x,y) & = \cases  0 & \hbox{ if } x \notin iD, \ y \neq x^* \\
                   1 & \hbox{ if } x \notin iD, \ y = x^* \\
                   p(x,y) & \hbox{ if } x \in iD \endcases.
\endaligned\tag{2.17}$$
The transition operator corresponding to the process $X_n$ killed upon
leaving $iD,$ denoted $T_D,$ is the map $T_D: C^0(iD) \to C^0(iD)$
defined by 
$$\align
T_Df(x) & = \sum_{y\in iD} q(x,y) f(y). \tag{2.18}
\endalign$$
\endproclaim

\proclaim{Proposition 2.5} Let $(G,\o,W)$ be a graph with 
geometry and suppose that $D$ is a domain of $G$ with nonempty
boundary.  Let $T_D$ be the transition operator for the natural random
walk killed on exiting $iD.$  Let $\pair{f}{g}_w$ be the inner product
defined by 
$$\align
\pair{f}{g}_{w} & = \sum_{x \in V} f(x) g(x) W_V(x) w_V(x)\tag{2.19}
\endalign$$
Then $T_D$ is self-adjoint with respect to the inner product on
$C^0(iD)$ defined by \rom{(2.19)} and the eigenvalues of $T_D$
all have magnitude less than $1.$  
\endproclaim

\demo{Proof} Let $f,g \in C^0(iD).$  Then
$$\align
\pair{T_Df}{g}_{w} & = \sum_{x\in iD} \sum_{y\in iD} f(y) q(x,y)
g(x) W_V(x)w_V(x) \\
 & = \sum_{x\in iD} \sum_{y\in iD} f(y) \frac{W_E(x,y)}{w_V(x)}g(x)
W_V(x)w_V(x) \\  
 & = \sum_{x\in iD} \sum_{y\in iD} f(y) \frac{W_E(y,x)}{w_V(y)}g(x)
W_V(y) w_V(y) = \pair{f}{T_Dg}_{w} 
\endalign$$
proving that $T_D$ is self-adjoint.  Let $q_n$ be defined recursively
by 
$$\align
q_1(x,y) & = q(x,y) \\
q_n(x,y) & = \sum_{z\in V} q(x,z) q_{n-1}(z,y).
\endalign$$
Since $P^x(X_n \in iD \ \ \forall n) = 0,$ we have that $q_n(x,y) \to
0$ for all $x,y \in D.$  From this we 
conclude that
$$\align 
& \lim_{n \to \infty} T^n_Df(x) = 0
\endalign$$
for every $x\in iD$ and every $f \in C^0(iD).$ This proves the second
claim of the proposition.
\enddemo

\proclaim{Corollary 2.1} Let $(G,\o,W)$ be a graph with 
geometry and suppose that $D$ is a domain of $G$ with nonempty
boundary.  Then the interior Laplace operator, $\did:C^0(iD) \to
C^0(iD)$ is invertible.
\endproclaim

\demo{Proof}  From Proposition 2.3 and Proposition 2.4 we have that 
$$\align
\did  & = \Omega_w(I- T_D). \tag{2.20}
\endalign$$
By Proposition 2.5, the eigenvalues of $I - T_D$ are contained in the
interval $(0,2).$  Since $\Omega_w$ is invertible, the corollary
follows. 
\enddemo

\proclaim{Proposition 2.6} Let $(G,\o,W)$ be a graph with 
geometry and suppose that $D$ is a domain of $G$ with nonempty
boundary.  Let $T_D$ be the transition operator for the natural random
walk killed on  exiting $iD$ and let $\eta$ be the first exit index
for the natural random walk on $G.$  Then for every positive integer
$l,$ 
$$\align
P^x(\eta = l) & = \left[ T_D^{l-1}(I - T_D)\right] (\chone_{iD})(x)
\tag{2.21} 
\endalign$$
where $P^x$ is the measure charging paths beginning at $x$ and
$\chone_{iD}$ is the indicator function of $iD.$
\endproclaim

\demo{Proof}  Note that $\chone_{iD} = \sum_{y \in iD} \delta_y$ where
$\delta_y$ is given by \rom{(2.6)}.  Hence,
$$\align
T_D\chone_{iD} (x) & = \sum_{z\in iD} \sum_{y\in iD} q(x,z)\delta_y(z) \\
 & = \sum_{y\in iD} q(x,y) 
\endalign$$
where, as above, $q(x,y)$ is the transition probability for the
natural random walk killed on leaving $iD.$  Since $\sum_{y\in iD}
q(x,y) $ is the probability that starting at $x$ we transition to a
point in $iD,$ we conclude that  
$$\align
[I - T_D] (\chone_{iD})(x) & = P^x(\eta = 1). \tag{2.22}
\endalign$$
The proposition now follows from \rom{(2.22)} and repeated use of the
fact that 
$$\align
P^x(\eta = 2) & = \sum_{y \in iD} q(x,y) P^y(\eta = 1).
\endalign$$
\enddemo
\proclaim{Definition 2.6} Suppose that $(G,\o,W)$ is a graph with
geometry and suppose that $D$ is a domain of $G$ with nonempty
boundary.  Let $T_D$ be the transition operator for the natural random
walk killed on exiting $iD.$  The Green operator of $iD$ is the
operator defined by 
$$\align
G & = \sum_{n=0}^\infty T_D^n. \tag{2.23}
\endalign$$
\endproclaim
Note that by Proposition 2.5 the sum in \rom{(2.23)} converges and $G$ is
well defined.  Indeed, we note that the Green operator gives the
inverse of the Laplace operator $\did:$  For every $f \in C^0(iD)$ and
every $x\in iD,$
$$\align
\did G \Omega_w^{-1} f(x) & = f(x). \tag{2.24}
\endalign$$

\heading
3:  Discrete Boundary Value Problems
\endheading

In this section we develop the machinery required to prove Theorem
1.1.  Throughout this section, $(G,\o,W)$ will be a graph with
geometry (cf Definition 2.1) and $D$ will be a domain of $G$ with
nonempty boundary (cf Definition 2.2).   

Let $\did$ be the interior Laplace operator
associated to $D$ (cf \rom{(2.10)}).  Define a sequence of functions
$f_k \in C^0(iD)$ recursively as follows:  Set 
$$\align
f_0(x) & = 1 \hbox{ for all } x\in iD \tag{3.1}
\endalign$$
and having defined $f_j, \ 1\leq j <k,$ let $f_k: iD \to \R$ be the unique
solution of 
$$\align
\did f_k + kf_{k-1} + \sum_{j=2}^{k} \choose{k}{j}
\left(-\frac{1}{w_V}\right)^{j-1}f_{k-j} & = 0 \tag{3.2} 
\endalign$$
where $w_V$ is the natural auxillary weighting for the graph
with geometry $(G,\o,W)$ (cf \rom{(2.2)}). 

The functions $f_k$ are closely related to the exit time moments (from
$D$) of the natural random walk associated to $(G, \o, W):$

\proclaim{Theorem 3.1} Let $D$ be a domain and, for $k$ a nonnegative
integer, let $f_k$ be defined by \rom{(3.1)} and \rom{(3.2)} above.
Let $X_n$ be the natural random walk on the underlying graph and
suppose that $P^x$ is the associated measure charging paths beginning
at $x.$  Let $\tau$ be the expected exit time (from $D$) of the
natural random walk as defined by \rom{(2.16)}.  Then   
$$\align
f_k(x) & = E^x[\tau^k] \tag{3.3}
\endalign$$
where $E^x$ is the expectation corresponding to the measure $P^x.$
\endproclaim

\demo{Proof}  From Proposition 2.3 and the definition of the
transition probabilities we have
$$\align
-\did E^x [\tau^k] & =  \left[\sum_{y\in iD} W_E(x,y) E^y[\tau^k] -
w_V(x) E^x[\tau^k] \right] \\
   & =  w_V(x) \left[\sum_{y\in iD} p(x,y) E^y[\tau^k] -
E^x[\tau^k] \right] . \tag{3.4}
\endalign$$
But  
$$\align
 E^y[\tau^k] & = \sum_{l=0}^\infty \left(\sum_{n=0}^{l-1}
\frac{1}{w_V(X_n)}\right)^k P^y(\eta = l)
\endalign$$
where $X_n$ is the natural random walk, $P^y$ is the measure charging
paths beginning at $y,$ and $\eta$ is as defined in \rom{(2.15)}.
Hence, 
$$\align
\sum_{y\in iD} p(x,y) E^y[\tau^k] & = \sum_{l=0}^\infty
\sum_{y\in iD} \left(\sum_{n=0}^{l-1} \frac{1}{w_V(X_n^y)}\right)^k
p(x,y)P^y(\eta = l) 
\endalign$$
where the $''y''$ superscript emphasizes that the paths contributing to
the sum begin at $y.$  Noting that for $l\geq 1,$
$$\align
P^x(\eta = l+1) & = \sum_{y\in iD} p(x,y) P^y(\eta = l),
\endalign$$
we have
$$\align
\sum_{y\in iD} p(x,y) E^y[\tau^k] & = \sum_{l=1}^\infty
\left(\left(\sum_{n=0}^{l} \frac{1}{w_V(X_n^x)}\right)
-\frac{1}{w_V(x)}\right)^k  P^x(\eta = l+1). \tag{3.5}
\endalign$$
From \rom{(3.4)} and \rom{(3.5)} we conclude that
$$\align
\did E^x[\tau^k] & =  w_V(x)\left[ E^x[(\tau -(w_V(x))^{-1})^k] - 
E^x[\tau^k] \right]. \tag{3.6}
\endalign$$
Theorem 3.1 follows immediately from \rom{(3.6)}.
\enddemo

The relationship between the solutions of the Poisson problems defined
by (\rom{3.1}) - \rom{(3.2)} and the exit time moments for the natural
random walk on the underlying graph are closely related to the
corresponding relationship between the behavior of Brownian motion in
a domain in a Riemannian manifold $M$ and the solutions of a hierarchy
of Poisson problems on the domain.  More precisely, suppose that 
$D\subset M$ is a smoothly bounded open domain with compact closure.
Let $\tau $ be the first exit time of classical Brownian motion from
$D.$  Suppose that $\Delta$ is the Laplace operator on $M$ and define
a sequence of functions, $u_k,$ recursively as follows:  Let
$$\aligned
u_0(x) & = 1 \hbox{ for all } x\in D \\
u_0(x) & = 0 \hbox{ for all } x\in \partial D
\endaligned\tag{3.7}$$
and, having defined $u_j$ for $0 \leq j <k,$ let $u_k$ be the solution
to the Poisson problem 
$$\aligned
\frac12 \Delta u_k + ku_{k-1} & = 0 \hbox{ for all } x\in D \\
u_k(x) & = 0 \hbox{ for all } x\in \partial D.
\endaligned\tag{3.8}$$
Then (cf \cite{H}, \cite{M})
$$\align
u_k(x) & = E^x[\tau^k] \tag{3.9}
\endalign$$
where $E^x$ is expectation with respect to the measure weighting
paths beginning at $x.$  

Given \rom{(3.9)}, the sum appearing in \rom{(3.2)} can be interpreted
as the correction necessary to insure that exit time moments for the
natural walk can be expressed as solutions to discrete Poisson
problems on the corresponding domain.  Note that this correction term
vanishes for the first moment.  

As mentioned in the introduction (cf \rom{(1.9)} and \rom{(1.10)}),
the quantities 
$$\align
B_{1,k}(D) & = \int_D E^x[\tau^k] dx \\
  & = \int_D u_k(x) dx
\endalign$$
are invariants of the underlying metric which can be computed via
variational principles.  In the case $k=1$ the Poisson problems which
compute the expected exit time for both the discrete case (cf
\rom{(3.2)}) and the smooth case (cf \rom{(3.8)}) have the same form,
and we expect to be able to formulate and prove the discrete analog of
\rom{(1.11)}.  When $k \neq 1,$ the structure of the corresponding
Poisson problems differ and we expect combinatorial corrections terms
to arise.  We proceed to determine the precise form of the
correction. 

Let $D$ be a domain in a graph with geometry and
define a sequence of functions recursively as follows:  Let 
$$\align
g_0 & = 1 \hbox{ for all } x\in iD \tag{3.10}
\endalign$$
and having defined $g_j$ for $0 \leq j <k,$ let $g_k$ be the unique
solution to the Poisson problem
$$\align
\did g_k + kg_{k-1} & = 0 \hbox{ for all } x\in iD. \tag{3.11}
\endalign$$

\proclaim{Theorem 3.2} Suppose that $D$ is a domain with nonempty
boundary and that $ \Delta_{iD}$ is the interior Laplace operator  
associated to $D.$  For $k$ a nonnegative integer, let $g_k$ be defined
by \rom{(3.10)}-\rom{(3.11)}.  Then, with the pairings defined by
\rom{(2.3)}, 
$$\align
\pair{g_k}{1}_{iD} & = k! \sup_{g\in C^0(iD)\setminus \{0\}}
\frac{\pair{g}{1}_{iD}^2}{|\pair{g}{\did^kg}_{iD}|}.
\endalign$$
\endproclaim

\demo{Proof}  The proof follows the argument given in \cite{KM}.  For
$g \in C^0(iD) \setminus\{0\},$  consider the quotient   
$$\align 
Q_k(g)  & =  \frac{ \pair{1}{g}_{iD}^2}{|\pair{g}{\did^kg}_{iD}|}.
\endalign$$
Note that $Q_k$ is invariant under scaling: for all nonzero scalars
$z,$ $$Q_k(zg) = Q_k(g).\tag{3.12}$$ 
Suppose $k = 2m.$  Since $\did$ is self-adjoint, we have
$$\align 
Q_k(g)  & =  \frac{ \pair{1}{g}_{iD}^2}{\|\did^mg\|_{iD}^2}.
\endalign$$
Since $\did$ is invertible we can write $g = \did^{-m}h$ for some $h \in
C^0(iD).$ Since $\did$ is self-adjoint and $Q_k$ is scale invariant
\rom{(3.12)}, we conclude that    
$$\align 
\sup_{g \in C^0(iD) \setminus \{0\}} Q_k(g)  & =  \sup_{h \in C^0(iD), \
\|h\|^2_{iD} = 1}  \pair{\did^{-m}(1)}{h}_{iD}^2. \tag{3.13}
\endalign$$ 
To maximize the inner product appearing on the right hand side of
\rom{(3.13)}, we take the unit vector in the direction of $\did^{-m}(1).$
That is, we set $h= \frac{\did^{-m}(1)}{\|\did^{-m}(1)\|},$ $g =
\frac{\did^{-2m}(1)}{\|\did^{-m}(1)\|},$ and use scale invariance to write
$$\align 
\sup_{g \in C^0(iD) \setminus \{0\}} Q_k(g)  & =
Q_k\left(\frac{\did^{-2m}(1)}{\|\did^{-m}(1)\|_{iD}}\right) \\
 & = \frac{\pair{1}{\did^{-2m}(1)}_{iD }^2}{|\pair{\did^{-2m}(1)
}{\did^{2m}\did^{-2m}(1)}_{iD}|} . \tag{3.14}   
\endalign$$
From the recursion \rom{(3.10)}-\rom{(3.11)} defining the $g_k$
we have 
$$\align
\did^{-k}(1) & = \frac{(-1)^k}{k!} g_k. \tag{3.15} 
\endalign$$
Using \rom{(3.15)}, we obtain
$$\align 
\sup_{g \in C^0(iD) \setminus \{0\}} Q_k(g)  & =
\frac{1}{(2m)!} \pair{1}{g_{2m}}_{iD }
\endalign$$
which concludes the proof of the theorem when $k$ is even.  The
argument for odd $k$ is similar:  Suppose that $k = 2m+1$ and write 
$$\align
Q_k(g)  & =  \frac{ \pair{1}{g}_{iD}^2}{|\pair{\did^mg}{\did \did^mg}_{iD}|}.
\endalign$$
Since $\did = d^*_W d,$ we have
$$\align
Q_k(g)  & =  \frac{ \pair{1}{g}_{iD}^2}{|\pair{d\did^mg}{d\did^mg}_{iD}|}
\tag{3.16} 
\endalign$$
where the pairing in the denominator is understood to be on edge
functions.  As above, set $g = \did^{-m}h$ and note that
$$\align 
\sup_{g \in C^0(iD)\setminus \{0\}} Q_k(g)  & =  \sup_{h \in C^0(iD), \
\|dh\|^2_{iD} = 1}  \pair{d\did^{-m-1}(1)}{dh}_{iD}^2. \tag{3.17}
\endalign$$ 
where once again the pairing is for edge functions.  As before, to
maximize the inner product appearing on the right hand side of 
\rom{(3.17)}, we take the unit vector in the direction of $d\did^{-m-1}(1).$
That is, we set $dh= \frac{d\did^{-m-1}(1)}{\|d\did^{-m-1}(1)\|}$ and use scale
invariance to write 
$$\align 
\sup_{g \in C^0(iD) \setminus \{0\}} Q_k(g)  & =
\frac{\pair{d\did^{-m-1}(1)}{d\did^{-m-1}(1)}_{iD
}^2}{|\pair{d\did^{-m-1}(1)}{d\did^{-m-1}(1)}_{iD}|} . \tag{3.18}    
\endalign$$
Integrating by parts and using \rom{(3.18)}, we obtain
$$\align 
\sup_{g \in C^0(iD) \setminus \{0\}} Q_k(g)  & = \frac{1}{(2m+1)!}
\pair{1}{g_{2m+1}}_{iD } 
\endalign$$
which concludes the proof of the theorem when $k$ is odd.  
\enddemo

\heading
4:  Proof of Theorem 1.1
\endheading

Throughout this section $(G,\o,W)$ will be a graph with natural
geometry and $D\subset G$ will be a domain with nonempty boundary
which is $\alpha$-weight regular (cf Definition \rom{(2.5)}).  

As above, let $X_n$ be the natural random walk on $G,$ let $\eta$ be
the first exit index of $X_n$ from $D$ (cf \rom{(2.15)}), and $\tau$
the first exit time of $X_n$ from $D$ (cf \rom{(2.16)}).  As a
consequence of $\alpha$-weight regularity, we have $$\eta = \alpha \tau.$$ 
To prove Theorem \rom{1.1} we will express the functions
$E^x[\tau^k]$ in terms of the functions $g_k(x)$ defined by
\rom{(3.10)}-\rom{(3.11)}.  With this in mind, consider the power series 
$$\align
P_k(x) & = x\sum_{l=1}^\infty l^k (1-x)^{l-1} \tag{4.1}
\endalign$$
which converges absolutely for $x \in (0,2).$  There is a
simple recursion for $P_k(x):$
$$\aligned
P_0(x) & = 1\\
P_{k+1}(x) & = \frac{1}{x} P_k(x) +(x-1) \ddx P_k(x)
\endaligned\tag{4.2}$$
The following definition will provide convenient notation:

\proclaim{Definition 4.1} Let $n$ and $k$ be positive integers.  The
Stirling numbers of the second kind are defined by 
$$\align
\sone{k}{n} & = \frac{(-1)^n}{n!} \sum_{m=1}^n \choose{n}{m} (-1)^m
m^k. \tag{4.3}
\endalign$$
The Stirling numbers of the first kind are defined by 
$$\align
\stwo{k}{n} & =  \sum_{m=0}^{k-n} (-1)^m \choose{k-1+m}{k-n+m}
\choose{2k -n}{k-n-m} \sone{k-n-m}{m}. \tag{4.4}
\endalign$$
\endproclaim

Recall that the Stirling numbers of the first kind count the number of
permutations of $k$ symbols which contain exactly $n$ symbols while
the Stirling numbers of the second kind represent the number of ways
of partitioning $k$ objects into $n$ blocks.  These numbers play a
role in combinatorics and statistical mechanics (cf \cite{Z}).

\proclaim{Lemma 4.1} For $k \geq 1,$ 
$$\align
P_k(x) & = (-1)^k \sum_{n=1}^k \sone{k}{n} n!
\left(-\frac{1}{x}\right)^n \tag{4.5}
\endalign$$
where $\sone{k}{n}$ is the Stirling number of the second kind.
\endproclaim

\demo{Proof} We check the case $k=1$ by starting with \rom{(4.1)},
summing the geometric series, and differentiating.  The general case
follows by induction using the recursion \rom{(4.2)}.
\enddemo

\proclaim{Lemma 4.2} Let $D$ be an $\alpha$-weight regular domain.
Then, with notation as above,
$$\align
E^x[\tau^k] & =  \sum_{n=1}^k \sone{k}{n} (-\alpha)^{n-k}
n! \did^{-n}(\chone_{iD})(x). \tag{4.6} 
\endalign$$
\endproclaim

\demo{Proof}  From the definition we have 
$$\align
E^x[\tau^k] & = \sum_{l=1}^\infty \left(\sum_{n=0}^{l-1}
\frac{1}{w_V(X_n)} \right)^k
P^x(\eta=l) \\
 & = \alpha^{-k} \sum_{l=1}^\infty l^kP^x(\eta = l).
\endalign$$
Using Proposition \rom{2.6} (cf \rom{(2.21)}), 
$$\align
E^x[\tau^k] & = \alpha^{-k} \sum_{l=1}^\infty l^k
[T_D^{l-1}(I-T_D)](\chone_{iD})(x). \tag{4.7} 
\endalign$$
By \rom{(2.20)} and $\alpha$-weight regularity, $$T_D = I -
\alpha^{-1} \did. \tag{4.8}$$
Using \rom{(4.7)}, \rom{(4.8)}, and \rom{(4.1)}, we obtain
$$\align
E^x[\tau^k] & = \alpha^{-k} \sum_{l=1}^\infty l^k
[(I - \alpha^{-1} \did)^{l-1}\alpha^{-1}\did (\chone_{iD})(x) \\
 & =  \alpha^{-k} P_k(\alpha^{-1}\did) (\chone_{iD})(x). \tag{4.9} 
\endalign$$
Using \rom{(4.9)} and Lemma \rom{4.1} we obtain 
$$\align
E^x[\tau^k] & = \alpha^{-k}(-1)^k \sum_{n=1}^k \sone{k}{n} n!
\left(-\frac{1}{\alpha^{-1}\did} \right)^n (\chone_{iD})(x) \tag{4.10} 
\endalign$$
from which the lemma follows.
\enddemo

\proclaim{Theorem 4.1}  Let $D$ be an $\alpha$-weight regular domain with
nonempty boundary.  Let $\tau$ be the first exit time of the natural
random walk from the domain $D$ and let $g_k$ be the solutions of the
Poisson problems defined recursively by \rom{(3.10)}-\rom{(3.11)}. Then
$$\align
E^x[\tau^k] & = g_k + \sum_{n=1}^{k-1} \sone{k}{n} (-\alpha)^{n-k} g_n
\tag{4.11} 
\endalign$$
where $\sone{k}{n}$ are the Stirling numbers of the second kind
\rom{(4.3)} 
\endproclaim

\demo{Proof} The recursion \rom{(3.10)}-\rom{(3.11)} defining $g_k$
gives 
$$\align
g_n(x) & = (-1)^n n! \did^{-n}(\chone_{iD})(x). \tag{4.12}
\endalign$$
Noting that $\sone{k}{k} = 1$ and using Lemma 4.2, the proof of the
theorem is complete.
\enddemo

\demo{Proof of Theorem 1.1}  Pairing both sides of \rom{(4.11)} gives
\rom{(1.7)}.  To see that \rom{(1.8)} holds, note that the Stirling
numbers satisfy the following:  If we define a polynomial $x^{(k)}$ by 
$$\align
x^{(k)} & = x(x-1)(x-2) \cdot (x-k+1)
\endalign$$
then (cf \cite{Z}),
$$\align
x^{(k)} & = \sum_{n= 1}^k \stwo{k}{n} x^n \tag{4.13}
\endalign$$
and 
$$\align
x^{k} & = \sum_{n= 1}^k \sone{k}{n} x^{(n)}. \tag{4.14}
\endalign$$
Thus, \rom{(1.8)} follows from \rom{(1.7)}, \rom{(4.13)} and
\rom{(4.14)}, which completes the proof of Theorem 1.1.
\enddemo

\heading
5:  Spectra and Graph Geometry
\endheading

In this section we give proofs for Theorem 1.2 and Theorem 1.3.  We
begin with a definition:  

\proclaim{Definition 5.1} Let $(G,\o,W)$ be a graph with geometry and
suppose that $D$ is a domain in $G$ with nonempty boundary.  Suppose
that $\spec(D)$ is the spectrum of the interior Laplace operator on $D$ and
suppose that  $\lambda_j \in \spec(D).$  Let $E_{\lambda_j}$ be the
eigenspace corresponding to $\lambda_j$ and let $\{\phi_{\lambda_j,l}
: 1\leq l \leq \text{\rm dim}(E_{\lambda_j})\}$ be an orthonormal basis of
$E_{\lambda_j}.$  Let 
$$\align
\hat{a}_j^2 & = \sum_{l=1}^{\text{\rm dim}(E_{\lambda_j})}
\pair{\phi_{j,l}}{\chone_{iD}}
\endalign$$
be the projection of $E_{\lambda_j}$ on the space of constant
functions.  We call the sequence 
$$\align
a^2 & = (\hat{a}_1^2, \hat{a}_2^2, \dots, \hat{a}_j^2,\dots) \tag{5.1}
\endalign$$
a spectral partition of the volume of $D.$ 
\endproclaim

We note that a spectral partition of volume is independent of the
choice of basis for each eigenspace. 

The following will be used on a number of occasions: 

\proclaim{Lemma 5.1}Suppose that $(G,\o,W)$ is a graph with geometry
and that $D$ is a domain of $G$ with nonempty boundary.  Let $g_k$ be
defined recursively as solutions to Poisson problems by
\rom{(3.10)}-\rom{(3.11)}.  Suppose that $\spec(D) $ is
the spectrum of the interior Laplace operator associated to $D$
and that $\phi_j$ is a normalized eigenfunction corresponding to the
eigenvalue $\lambda_j$ (where we 
take an orthonormal basis of eigenfunctions when the multiplicity is
greater than one).  Let $$a_j = \pair{\chone_{iD}}{\phi_j}_V\tag{5.2}$$
where $\chone_{iD}$ is the indicator function of $iD.$  Then  
$$\align
\pair{g_k}{\phi_j}_V  & = \frac{(-1)^kk!}{\lambda_j^k}a_j. \tag{5.3}
\endalign$$
\endproclaim

\demo{Proof} With notation as above, we have $$\chone_{iD} = \sum_j a_j
\phi_j.\tag{5.4}$$  Using that $g_k$ is defined by the recursion
\rom{(3.10)}-\rom{(3.11)} and that $\did$ is self-adjoint, we obtain 
$$\align 
\pair{g_k}{\phi_j}_V & = \lambda_j^{-1}\pair{g_k}{L\phi_j}_V \\
 & = \frac{-k}{\lambda_j}\pair{g_{k-1}}{\phi_j}_V
\endalign$$
Continuing inductively and using \rom{(5.4)} we obtain \rom{(5.3)}. 
\enddemo

\proclaim{Corollary 5.1} Let $(G,\o,W)$ be a graph with geometry and
suppose that $D$ is a domain in $G$ with nonempty boundary.  Let
$\spec(D)$ be the spectrum of the interior Laplace operator on $D$ and
suppose that $\spec^*(D)$ is given by \rom{(1.12)}.  Let $\mspec(D)$
and $\pspec(D)$ be given by \rom{(1.5)} and \rom{(1.6)}, respectively.
Then $\pspec(D)$ is determined by $\spec^*(D)$ and the spectral
partition of the volume.  When $D$ is $\alpha$-weight regular, the 
same claims hold with $\pspec(D)$ replaced by $\mspec(D).$ 
\endproclaim

\demo{Proof}For each distinct eigenvalue, fix an orthonormal basis of
the corresponding eigenspace.  Let $a^2 = (\hat{a}_1^2, \hat{a}_2^2,
\dots, \hat{a}_j^2,\dots) $ be a spectral partition of volume as
defined in \rom{(5.1)}.  From Lemma 5.1 and Theorem 3.2, we have 
$$\align
A_{2,k} & = \sum_j \left( -\frac{1}{\lambda_j}\right)^k
k!\hat{a}_j^2. \tag{5.5} 
\endalign$$
where the sum is over elements of $\spec^*(D).$  The claims of the
corollary concerning elements of $\pspec(D)$ follow 
from \rom{(5.5)}.  By Theorem 1.1, the same is true for $\mspec(D)$
when $D$ is $\alpha$-weight regular.  
\enddemo

In what follows, we construct a partial converse for Corollary
\rom{5.1}.   

\demo{Proof of Theorem 1.2} Let $\{\mu_j\}$ be positive real numbers.
The {\it Problem of Moments} as formulated by Stieltjes asks:  For
which sequences $\{\mu_j\}$ is it possible to find a bounded
nondecreasing function $\psi: [0,\infty) \to \R$ such that 
$$\align
\mu_j & = \int_0^\infty x^j d\psi(x) ?  \tag{5.6}
\endalign$$
The problem arises in a variety of contexts (probability and
statistics, orthogonal polynomials, mechanics, etc) and has an
extensive associated literature which begins with a series of papers
by Tchebycheff starting in 1855 (cf \cite{SH} for background on the
moment problem).  In 1894-95, Stieltjes gave a treatment of the
problem (in which he developed the notion of the Stieltjes integral)
and clarified under what conditions the problem of moments admits a
solution.  We recall the relevant ideas and notation.

For $n$ a nonnegative integer, let $M_{0,n}$ be the $n \times
n$ symmetric square matrix defined by 
$$\align
(M_{0,n})_{ij} & = \mu_{i+j}, \ \ \ 0 \leq i,j \leq n-1. \tag{5.7}
\endalign$$
Let 
$$\align
(M_{1,n})_{ij} & = \mu_{i+j+1}, \ \ \ 0 \leq i,j \leq n-1. \tag{5.8}
\endalign$$
We have the following:
\proclaim{Theorem}(cf \cite{SH}) Let $\{\mu_j\}$ be positive real
numbers.  A necessary condition for the existence of a solution of the
Stieltjes moment problem 
$$\align
\mu_j & = \int_0^\infty x^j d\psi(x) 
\endalign$$
is that 
$$\align
\deter(M_{0,n}) \geq 0 , \ \ \  & \deter(M_{1,n}) \geq 0  \hbox{ for all }
n.
\endalign$$
In order that the spectrum of the problem not reduce to a finite set
it is necessary and sufficient that 
$$\align
\deter(M_{0,n})> 0 , \ \ \  & \deter(M_{1,n}) > 0  \hbox{ for all }
n.
\endalign$$
In order that there exists a solution whose spectrum is exactly $k+1$
points distinct from $0,$ it is necessary and sufficient that 
$$\aligned
\deter(M_{0,n})> 0 , \ \ \  & \deter(M_{1,n}) > 0  \hbox{ for all }
n, \ 0 \leq n \leq k \\
\deter(M_{0,n})= 0 , \ \ \  & \deter(M_{1,n}) = 0  \hbox{ for all } n > k.
\endaligned\tag{5.9}$$
In the last case, $\psi $ is uniquely determined by the sequence
$\{\mu_j\}.$  
\endproclaim
For $N$ fixed, let ${\Cal D}_N$ be as defined in \rom{(1.13)}.  Fix
$D \in {\Cal D}_N$ and choose $$\mu_n =\frac{A_{2,n}}{n!}.$$
Let $\hat{a}_j^2$ be the elements of a spectral partition of volume as
defined in \rom{(5.1)}.  Then, from \rom{(5.5)}, we have 
$$\align
\mu_n & = \sum_{j=1}^N \left(-\frac{1}{\lambda_j}\right)^n
\hat{a}_j^2\tag{5.10}  
\endalign$$
where the sum runs over $j$ for which $\lambda_j \in \spec^*(D).$
Define $\psi(x)$ by    
$$\align
\psi(x) & = \sum_{j=1}^N \hat{a}_j^2 \chone_{[-\frac{1}{\lambda_j},
\infty)}(x)\tag{5.11}
\endalign$$
Clearly, $\psi(x)$ is the unique solution to the moment problem and
hence the conditions given in \rom{(5.9)} hold where $N$ is the number
of elements in the set $\spec^*(D).$  Since $\pspec(D)$ determines the
set $\{\mu_j\},$ this implies that the Poisson spectrum determines
$\spec^*(D),$ as claimed in Theorem 1.2.  By Theorem 1.1, when $D$ is
$\alpha$-weight regular, the same is true for $\mspec(D).$    

To understand the precise relationship between the Poisson spectrum
and $\spec^*(D),$ Let $X $ be the vector $(X_1,X_2, \dots , X_N)$ and
consider the Vandermonde matrix 
$$\aligned
V_X & = \left( \matrix
                   1 & 1 & \dots & 1  \\
                   X_1 & X_2& \dots & X_N \\
                   \vdots & \vdots & \vdots & \vdots \\
                   X_1^{N-1} & X_2^{N-1}& \dots & X_N^{N-1}
               \endmatrix
          \right).
\endaligned\tag{5.12}$$
Note that when $X_i$ are distinct, $V_X$ is invertible.  Let
$\mu_n$ be defined as in \rom{(5.10)} and set $$\mu =
(\mu_0,\mu_1,\dots,\mu_{N-1}).\tag{5.13}$$Let $\{\lambda_j\} = \spec^*(D)$
and write $$-\lambda^{-1} =
\left(-\frac{1}{\lambda_1},-\frac{1}{\lambda_2},\dots,-\frac{1}{\lambda_N}
\right).\tag{5.14}$$
Let $a^2$ be the spectral partition of volume defined as in
\rom{(5.1)}: $$a^2 = (\hat{a}_1^2,\hat{a}_2^2, \dots 
\hat{a}_N^2).\tag{5.15}$$
By \rom{(5.10)} we have  
$$\align
\mu & = V_{-\lambda^{-1}} a^2.\tag{5.16}
\endalign$$
Let $-\Lambda^{-1}$ be the diagonal matrix defined by
$-\lambda^{-1}.$  Let $$b_j = (\mu_j, u_{j+1}, \dots,
u_{j+N-1})^T$$ where the superscript denotes transpose.  Then
$b_n = V_{-\lambda^{-1}}(-\Lambda^{-1})^n a^2$ and, since the
$\lambda_j$ are distinct (this is the point of introducing the
spectral partition of volume), $V_{-\lambda^{-1}}$ is invertible.
Thus, we have   
$$\align
b_j & = V_{-\lambda^{-1}}(-\Lambda^{-1}) V_{-\lambda^{-1}}^{-1}
b_{j-1}. \tag{5.17} 
\endalign$$
Thinking of the vectors $b_n$ as columns of a matrix and using
\rom{(5.17)} repeatedly, we obtain
$$\align
M_{1,N} & = V_{-\lambda^{-1}} (-\Lambda^{-1}) V_{-\lambda^{-1}}^{-1}
M_{0,N} \tag{5.18}
\endalign$$
where $M_{0,N}$ and $M_{1,N}$ are given by \rom{(5.7)} and
\rom{(5.8)}, respectively.  Since $M_{0,N}$ is invertible we obtain 
$$\align
V_{-\lambda^{-1}} (-\Lambda^{-1}) V_{-\lambda^{-1}}^{-1} & =
M_{1,N}(M_{0,N})^{-1}.  
\tag{5.19} 
\endalign$$
The characteristic polynomial of $-\Lambda^{-1}$ is invariant under
conjugation by an invertible matrix.   Thus, the polynomial 
$$\align
P_{D} (x) & = \prod_{j} \left(- \frac{1}{\lambda_j} - x \right) \\
   & = \deter( M_{1,N} (M_{0,N})^{-1} - xI) \tag{5.20}
\endalign$$ 
has coefficients which are polynomial in the entries of
the matrix $ M_{1,N} (M_{0,N})^{-1}.$  In particular, we can conclude
that $\spec^*(D)$ is determined by the roots of a polynomial whose
coefficients are rational functions of the $\mu_j.$  We conclude that
the polynomial appearing in \rom{(5.20)} satisfies the claims of
Theorem \rom{1.2}, which concludes the proof of the theorem.  
\enddemo

\proclaim{Corollary 5.2}  Let $(G,\o,W)$ be a graph with geometry and
suppose that $D$ is a domain in $G$ with nonempty boundary.  Let $a^2$
be the spectral partition of volume for $D$ as defined in
\rom{(5.1)}.  Then $\pspec(D)$ determines $a^2.$ 
\endproclaim

\demo{Proof} By Theorem 1.2, $\pspec(D)$ determines $\spec^*(D).$
Hence, $\pspec(D)$ determines the Vandermonde matrix given by
\rom{(5.12)}, as well as its inverse.  The corollary now follows from
\rom{(5.16)}. 
\enddemo

\proclaim{Definition 5.2}Suppose that $(G,\o,W)$ is a graph with
geometry and that $D$ is a domain of $G.$  Suppose that $H(x,t)$ is a
solution to the following boundary value problem with initial
data: 
$$\aligned
\did H & = \partial_t H \hbox{ on } iD\times (0,\infty) \\
H(x,0) & = 1 \hbox{ on } iD \\
H(y,t) & = 0 \hbox{ on } \partial D \times (0,\infty)
\endaligned\tag{5.21}$$
The {\it heat content of} $D$ is the function $Q(t)$ defined by 
$$\align
Q(t) & = \pair{H}{\chone_{iD}}_V \tag{5.22}
\endalign$$
where $\chone_{iD}$ is the indicator function of $iD.$
\endproclaim
The heat content admits a power series expansion:
$$\align
Q(t) & \simeq \sum_{n=0}^\infty q_n t^n. \tag{5.23}
\endalign$$
The coefficients $q_n$ in the expansion \rom{(5.23)} are called
{\it the heat content asymptotics for the domain} $D.$  In the
remainder of this paper, we investigate the relationship between heat
content asymptotics and spectra.

\demo{Proof of Theorem 1.3} Let  $(G,\o,W)$ be a graph with
geometry and suppose that $D$ is a domain of $G$ with nonempty
boundary.  Suppose that $H(x,t)$ solves \rom{(5.21)}.  Then, with the
eigenvalues of $\did$ denoted by $\lambda_j,$ the corresponding
normalized eigenfunctions denoted by $\phi_j,$ and $a_j$ defined by
\rom{(5.2)}, we have  
$$\align 
H(x,t) & = \sum_j a_j\phi_j(x) e^{\lambda_j t} \tag{5.24}
\endalign$$
where the sum is over $j$ with $\lambda_j \in \spec(D).$  Thus,
$$\align
Q(t)  & = \sum_j \hat{a}_j^2 e^{\lambda_j t}. \tag{5.25}
\endalign$$
where $\hat{a}_j^2$ are elements of a spectral partition of volume and
the sum is over $\lambda_j \in \spec^*(D).$  Using \rom{(5.25)} and
power series expansions of $e^{\lambda_j t}$ we 
see that the heat content asymptotics are given by 
$$\aligned
 q_n & = \sum_j \hat{a}_j^2 \frac{\lambda_j^n}{n!} \\
\endaligned\tag{5.26}$$
where once again the sum is over $j$ with $\lambda_j \in \spec^*(D).$
The second claim of Theorem 1.3 now follows from \rom{(5.26)}, Theorem
1.2, and Corollary 5.2.

Given \rom{(5.3)}, and the similarity to \rom{(5.10)}, it is
instructive to follow the technique introduced in the proof of Theorem
1.2 to study the relationship of the heat content asymptotics to
$\spec^*(D).$  To this end, let $$\mu_n = (-1)^nn!q_n.$$Then $\mu_n$
is positive and we have
$$\align
\mu_n & = \sum_j (-\lambda_j)^n \hat{a}_j^2. \tag{5.27}
\endalign$$
Let $-\lambda = (-\lambda_1,-\lambda_2, \dots -\lambda_N),$ let $a^2$ be
defined as in \rom{(5.1)} and let $\mu$ be defined as in
\rom{(5.13)}.  Then we have $$\mu = V_{-\lambda} a^2$$ where $V_{-\lambda}$
is the Vandermonde matrix given by \rom{(5.12)}.  Let
$\Lambda_{-\lambda}$ be the diagonal matrix defined by $-\lambda.$
Setting $c_j = (\mu_j,\mu_{j+1},\dots, \mu_{j+N})^T,$ assuming that
$V_{-\lambda},$ is invertible and following our previous computation,
we have  
$$\align
c_j & = V_{-\lambda} \Lambda_{-\lambda} (V_{-\lambda})^{-1} c_{j-1}
\endalign$$
which leads immediately to 
$$\align
M_{1,N} & = V_\lambda \Lambda_\lambda V_{\lambda}^{-1} M_{0,N}
\endalign$$
where $M_{0,N}$ and $M_{1,N}$ are given by \rom{(5.6)} and
\rom{(5.7)}, respectively.  As before, we can construct an explicit
solution to the corresponding Stieltjes Moment Problem:
$$\align
\psi(x) &  = \sum_j \hat{a}_j^2 \chone_{[-\lambda_j, \infty)}(x).
\endalign$$ 
As before, when $M_{0,N}$ is invertible, we obtain 
$$\align
V_\lambda \Lambda_\lambda V_{\lambda}^{-1} & = M_{1,N}(M_{0,N})^{-1}
. \tag{5.28} 
\endalign$$
From \rom{(5.28)} we conclude that the heat content asymptotics
determine $\spec^*(D),$ proving the second claim of Theorem 1.3.  Our
computations allow us to conclude that, as in the statement of Theorem
1.2, for $N$ fixed, there are $N$ rational functions $f_i', \ 1 \leq i
\leq N, \ f_i':\R^{2N} \to \R$ such that for every $D \in {\Cal
D}_N,$ the roots of the polynomial $$P_N'(x) = x^N + \sum_{i=0}^{N-1}
f_i'(A_{2,0}(D), A_{2,1}(D), \dots, A_{2,2N-1}(D))x^i \tag{5.29} $$
give the elements of $\spec^*(D),$ which concludes the proof of
Theorem 1.3.  
\enddemo 

We note that the heat content asymptotics are determined by
$\spec^*(D)$ and the spectral partition of volume:

\proclaim{Corollary 5.3} Let $(G,\o,W)$ be a graph with geometry and
suppose that $D$ is a domain in $G$ with nonempty boundary.  Let
$\spec(D)$ be the spectrum of the interior Laplace operator on $D$ and
suppose that $\spec^*(D)$ is given by \rom{(1.12)}.  Let $\{q_n \}$ be
the heat content asymptotics of $D.$  Then $\{q_n\}$ is determined by
$\spec^*(D)$ and the spectral partition of the volume.  
\endproclaim

\demo{Proof}  This follows immediately from \rom{(5.26)}.
\enddemo

Finally, we note that the correspondence ``$\lambda \leftrightarrow
\frac{1}{\lambda}$'' interchanges the role of heat content asymptotics
and Poisson spectrum.  We can formalize this using a construction
closely related to that of spectral zeta functions:

\proclaim{Definition 5.4} Suppose that $(G,\o,W)$ is a graph with
geometry and that $D$ is a domain of $G$ with nonempty boundary.  Let
$\spec(D)$ be the Dirichlet spectrum associated to $D$ and let $a_j$
be defined as in \rom{(5.2)}.  Let $s$ be a complex variable.  The
weighted zeta function associated to $D$ is defined to be
$$\align
\zeta_{D}(s) & = \sum_j a_j^2 \left(- \frac{1}{\lambda_j}\right)^s.
\tag{5.30} 
\endalign$$
\endproclaim
We note that while $a_j$ depends on the choice of basis for each
eigenspace, the definition of $\zeta_D(s)$ does not.  

\proclaim{Corollary 5.4}  Suppose that $(G,\o,W)$ is a graph with 
geometry and that $D$ is a domain of $G$ with nonempty boundary.  Let
$\zeta_D(s)$ be the zeta function associated to $D$ as in
\rom{(5.30)}.  Then, with $n$ a positive integer, with $A_{2,n}$ as in
\rom{(1.6)}, and with $q_n$ as in \rom{(1.17)},
$$\align
\zeta_D(n) & = \frac{A_{2,n}}{n!} \tag{5.31} \\  
\zeta_D(-n) & =(-1)^n q_n n!. \tag{5.32}
\endalign$$
\endproclaim

\demo{Proof} The identity \rom{(5.31)} follows from \rom{(5.5)} and
the definition of the zeta function.  Similarly, the identity
\rom{(5.32)} follows from \rom{(5.26)}.
\enddemo

\enddocument